\documentclass[final]{siamltex}
\usepackage{a4} 
\usepackage{amsmath}
\usepackage{amssymb}
\usepackage{graphicx}
\usepackage{caption,subcaption}
\usepackage{multirow}
\usepackage{xstring}
\usepackage[utf8]{inputenc}
\usepackage{amstext,amsbsy,amsopn,eucal,enumerate}
\usepackage{tikz}
\usepackage{pgfplots}
\usepackage{tabularx}

\newcommand{\beq}{\begin{equation}}
\newcommand{\eeq}{\end{equation}}
\newcommand {\mat}      [1] {\left[\begin{array}{#1}}
\newcommand {\rix}          {\end{array}\right]}
\newcommand {\smat}      [1] {\left[\begin{smallmatrix}{#1}}
\newcommand {\srix}          {\end{smallmatrix}\right]}
\newcommand {\s}      [1] {\begin{smallmatrix}{#1}}
\newcommand {\se}          {\end{smallmatrix}}

\begin{document}

\title{A new type of singular perturbation approximation for stochastic bilinear systems}

\author{Martin Redmann\thanks{Weierstrass Institute for Applied Analysis and Stochastics, Mohrenstrasse 39, 10117 Berlin Germany  (Email: {\tt 
martin.redmann@wias-berlin.de})}. The author~gratefully acknowledge the support from the DFG through the research unit FOR2402.}

\maketitle

\begin{abstract}
Model order reduction (MOR) techniques are often used to reduce the order of spatially-discretized (stochastic) partial differential equations and hence reduce 
computational complexity. A particular class of MOR techniques is balancing related methods which rely on simultaneously diagonalizing the system Gramians. This has been extensively studied for deterministic linear systems.
The balancing procedure has already been extended to bilinear equations \cite{typeIBT}, an important subclass of nonlinear systems. The choice of Gramians in \cite{typeIBT} is referred to be the standard approach. 
In \cite{hartmann}, a balancing related MOR scheme for bilinear systems called singular perturbation approximation (SPA) has been described that relies on the standard choice of Gramians. However, no error bound for this method could be proved.
In this paper, we extend the setting used in \cite{hartmann} by considering a stochastic system with bilinear drift and linear diffusion term. Moreover, we propose a modified reduced order model and choose a different reachability 
Gramian. Based on this new approach, an $L^2$-error bound is proved for SPA which is the main result of this paper. This bound is new even for deterministic bilinear systems.
\end{abstract}

\begin{keywords}
 model order reduction, singular perturbation approximation, nonlinear stochastic systems, L\'evy process
\end{keywords}

\begin{AMS}
Primary: 93A15, 93C10, 93E03. Secondary: 15A24, 60J75. 
\end{AMS}


\section{Introduction}

Many phenomena in real life can be described by partial differential equations (PDEs). For an accurate mathematical modeling of these real world 
applications, it is often required to take random effects into account. Uncertainties in a PDE model can, for example, be represented by an 
additional noise term leading to stochastic PDEs (SPDEs) \cite{dapratozab,newspde,zabczyk,prevotroeckner}.\smallskip

It is often necessary to numerically approximate time-dependent SPDEs since analytic solutions do not exist in general. Discretizing in space can be 
considered as a first step. This can, for example, be done by spectral Galerkin \cite{galerkin, galerkinhaus, galerkinjentzen} or finite 
element methods \cite{MR1637047,MR2646102,Kruse14}. This usually leads to large-scale SDEs. 
Solving such complex SDE systems causes large computational cost. In this context, model order reduction (MOR) is used to save computational time by
replacing high dimensional systems by systems of low order in which the main information of the original system should be captured. 

\subsection{Literature review}

Balancing related MOR schemes were developed for deterministic linear systems first. Famous representatives of this class of methods are balanced truncation (BT) \cite{antoulas,moore,obiand} and 
singular perturbation approximation (SPA) \cite{fernic,spa}. \smallskip

BT was extended in \cite{bennerdamm, redmannbenner} and SPA was generalized in \cite{redSPA} to stochastic linear systems.
With this first extension, however, no $L^2$-error bound can be achieved \cite{bennerdammcruz, dammbennernewansatz}. 
Therefore, an alternative approach based on a different reachability Gramian was studied for stochastic linear systems leading to an $L^2$-error bound for BT \cite{dammbennernewansatz} and for SPA \cite{redmannspa2}.\smallskip

BT \cite{typeIBT,bennerdamm} and SPA \cite{hartmann} were also generalized to bilinear systems, which we refer to as the standard approach for these systems. Although bilinear terms are very weak nonlinearities, 
they can be seen as a bridge between linear and nonlinear systems. This is because many nonlinear systems can be represented by bilinear systems using a so-called Carleman linearization. Applications of these equations can
be found in various fields \cite{brunietal,Mohler,rugh}. The standard approach for bilinear system has the drawback that no $L^2$-error bound could be shown so far. A first error bound for the standard ansatz 
was recently proved in \cite{beckerhartmann}, where an output error bound in  $L^\infty$ was formulated for infinite dimensional bilinear systems. Based on the alternative choice of Gramians in \cite{dammbennernewansatz},
a new type of BT for bilinear systems was considered \cite{redmanntypeiibilinear} providing an $L^2$-error bound under the assumption of a possibly small bound on the controls.\smallskip

A more general setting extending both the stochastic linear and the deterministic bilinear case was investigated in \cite{redstochbil}. There, BT was studied and an $L^2$-error bound was proved overcoming the restriction 
of bounded controls in \cite{redmanntypeiibilinear}. In this paper, we consider SPA for the same setting as in \cite{redstochbil} in order to generalize the work in \cite{hartmann}. Moreover, we modify the reduced order 
model (ROM) in comparison to \cite{hartmann} and show an $L^2$-error bound which closes the gap in the theory in this context.\smallskip  

For further extensions of balancing related MOR techniques to other nonlinear systems, we refer to \cite{bennergoyal, Scherpen}.

\subsection{Setting and ROM}\label{settingstochstabgen}

Let every stochastic process appearing in this paper be defined on a filtered probability space $\left(\Omega, \mathcal F, \left(\mathcal F_t\right)_{t\geq 0}, \mathbb 
P\right)$\footnote{We assume that $\left(\mathcal F_t\right)_{t\geq 0}$ is 
right-continuous and $\mathcal F_0$ contains all sets $A$ with $\mathbb P(A)=0$.}. 
Suppose that $M=\left(M_1, \ldots, M_{v}\right)^T$ is an $\left(\mathcal F_t\right)_{t\geq 0}$-adapted and $\mathbb R^{v}$-valued mean zero L\'evy 
process with $\mathbb E \left\|M(t)\right\|^2_2=\mathbb E\left[ M^T(t)M(t)\right]<\infty$ for all $t\geq 0$. Moreover, we assume that for all $t, h\geq 0$ the random 
variable $M\left(t+h\right)-M\left(t\right)$ is independent of $\mathcal F_t$.\smallskip

We consider a large-scale stochastic control system with bilinear drift that can be interpreted as a spatially-discretized SPDE. We investigate the system
     \begin{subequations}\label{controlsystemoriginal}
\begin{align}
d x(t)&=[A x(t)+ Bu(t) + \sum_{k=1}^m N_k x(t) u_k(t)]dt+ \sum_{i=1}^v H_i x(t-) dM_i(t),\label{stateeq}\\ 
y(t)&= {C} x(t),\;\;\;t\geq 0.\label{originalobserveq}
\end{align}
\end{subequations}
We assume that $A, N_k, H_i\in \mathbb R^{n\times n}$ 
($k\in\left\{1, \ldots, m\right\}$ and $i\in\left\{1, \ldots, v\right\}$), $B\in \mathbb R^{n\times m}$ and $C\in \mathbb R^{p\times n}$. Moreover, we define $x(t-):=\lim_{s\uparrow t} x(s)$.
The control $u=\left(u_1, \ldots, u_{m}\right)^T$ is assumed to be deterministic and square integrable, i.e., 
     \begin{align*}
\left\|u\right\|_{L^2_T}^2:=\int_0^T \left\|u(t)\right\|_2^2 dt<\infty
\end{align*}                            
for every $T>0$. By \cite[Theorem 4.44]{zabczyk} there is a matrix $K=\left(k_{ij}\right)_{i, j=1, \ldots, v}$ such that  
$\mathbb E[M(t)M^T(t)]=K t$. $K$ is called covariance matrix of $M$. \smallskip 

In this paper, we study SPA to obtain a ROM. SPA is a balancing related method and relies on defining a reachability Gramian $P$ and an observability Gramian $Q$.
These matrices are selected, such that $P$ characterizes the states in (\ref{stateeq}) and $Q$ the states in (\ref{originalobserveq}) which barely contribute to the system dynamics, see \cite{redstochbil} for 
estimates on the reachability and observability energy. The estimates in \cite{redstochbil} are global, whereas the standard choice of Gramians leads to results being valid in a small neighborhood of 
zero only \cite{bennerdamm, graymesko}. 

In order to ensure the existence of these Gramians, throughout the paper it is assumed that 
\begin{align}\label{stochstab}
\lambda\left(A\otimes I+I\otimes A+\sum_{k=1}^m N_k\otimes N_k+\sum_{i, j=1}^v H_i\otimes H_j k_{ij}\right)\subset \mathbb C_-.
\end{align}
Here, $\lambda\left(\cdot\right)$ denotes the spectrum of a matrix. The reachability Gramian $P$ and the observability Gramian $Q$ are, according to \cite{redstochbil}, defined as the solutions to 
\begin{align}\label{newgram2}
 A^T P^{-1}+P^{-1}A+\sum_{k=1}^m N^T_k P^{-1} N_k + \sum_{i, j=1}^v H_i^T P^{-1} H_j k_{i j} &\leq -P^{-1}BB^T P^{-1},\\
 \label{gengenlyapobs}
A^T Q+Q A+\sum_{k=1}^m N_k^T Q N_k +\sum_{i, j=1}^v H_i^T Q H_j k_{ij} &\leq -C^T C,
                                       \end{align}
where the existence of a positive definite solution to (\ref{newgram2}) goes back to \cite{dammbennernewansatz, redmannspa2}.\smallskip

We approximate the large scale system (\ref{controlsystemoriginal}) by a system which has a much smaller state dimension $r\ll n$. This reduced order model (ROM) 
is supposed be chosen, such that the corresponding output $y_r$ is close to the original one, i.e., $y_r\approx y$ in some metric. 
In order to be able to remove both the unimportant states in (\ref{stateeq}) and (\ref{originalobserveq}) 
simultaneously, the first step of SPA is a state space transformation 
\begin{align*}
 (A, B, C, H_i, N_k)\mapsto (\tilde A, \tilde B, \tilde C, \tilde H_i, \tilde N_k):=(SAS^{-1}, SB, CS^{-1}, SH_iS^{-1}, SN_kS^{-1}),
\end{align*}
where $S=\Sigma^{-\tfrac{1}{2}} X^T L_Q^T $ and $S^{-1}=L_PY\Sigma^{-\tfrac{1}{2}}$. The ingredients of the balancing transformation are computed by the Cholesky factorizations $P=L_PL_P^T$, 
$Q=L_QL_Q^T$, and the singular value decomposition $X\Sigma Y^T=L_Q^TL_P$. This transformation does not change the output $y$ of the system, but it guarantees that the new Gramians are
diagonal and equal, i.e., $S P S^T=S^{-T}Q S^{-1}=\Sigma=\diag(\sigma_1,\ldots, \sigma_n)$ with $\sigma_1\geq \ldots \geq \sigma_n$ being the Hankel singular values (HSVs) of the system. \smallskip

We partition the balanced coefficients of (\ref{controlsystemoriginal}) as follows:
\begin{align}\label{partvonorimodel}
 \tilde A=\smat{A}_{11}&{A}_{12}\\ 
{A}_{21}&{A}_{22}\srix,\;\tilde B=\smat B_1 \\ B_2\srix,\; \tilde N_k=\smat{N}_{k, 11}&{N}_{k, 12}\\ 
{N}_{k, 21}&{N}_{k, 22}\srix,\;\tilde H_i=\smat{H}_{i, 11}&{H}_{i, 12}\\ 
{H}_{i, 21}&{H}_{i, 22}\srix,\;\tilde C= \smat C_1 & C_2\srix,
            \end{align}
where $A_{11}, N_{k, 11}, H_{i, 11}\in \mathbb R^{r\times r}$  ($k\in\left\{1, \ldots, m\right\}$ and $i\in\left\{1, \ldots, v\right\}$), $B_1\in \mathbb R^{r\times m}$ and $C_1\in \mathbb R^{p\times r}$ etc.
Furthermore, we partition the state variable $\tilde x$ of the balanced system and the diagonal matrix of HSVs \begin{align}\label{partitionfullmodel}
               \tilde x=\mat{c} x_1 \\ x_2\rix\text{ and }\Sigma=\mat{cc} \Sigma_1& \\ & \Sigma_2\rix,
                         \end{align}
where $x_1$ takes values in $\mathbb R^r$ ($x_2$ accordingly), $\Sigma_1$ is the diagonal matrix of large HSVs and $\Sigma_2$ contains the small ones.

Based on the balanced full model (\ref{controlsystemoriginal}) with matrices as in (\ref{partvonorimodel}), the ROM is obtained by neglecting the state variables $x_2$ corresponding to the small HSVs.
The ROM using SPA is obtained by setting $dx_2(t)=0$ and furthermore neglecting the diffusion and bilinear term in the equation related to $x_2$. The resulting 
algebraic constraint can be solved and leads to $x_2(t)=-A_{22}^{-1}(A_{21} x_1(t)+B_2 u(t))$. Inserting this expression into the equation for $x_1$ and into the output equation, the reduced system is 
\begin{subequations}\label{romstochstatebt}
\begin{align}\label{romstateeq}
             dx_r&=[\bar A x_r+\bar B u+\sum_{k=1}^m (\bar N_{k} x_r + \bar E_{k} u)u_k]dt+\sum_{i=1}^v (\bar H_{i} x_r+ \bar F_{i} u) dM_i,\\ 
    y_r(t)&=\bar Cx_r(t)+ \bar D u(t), \;\;\;t\geq 0,
            \end{align}
            \end{subequations}
with matrices defined by           
\begin{align*}
      &\bar A:=A_{11}- A_{12} A_{22}^{-1} A_{21},\;\;\;\bar B:=B_1-A_{12}A_{22}^{-1} 
B_2,\;\;\;\bar C:=C_1-C_2 A_{22}^{-1} A_{21},\\ &\bar D:=-C_2 A_{22}^{-1} B_2,\;\;\; \bar E_{k}:=-N_{k, 12} A_{22}^{-1} B_2,\;\;\; \bar F_{i}:=-H_{i, 12} A_{22}^{-1} B_2,\\
&\bar H_i:=H_{i, 11}-H_{i, 12} A_{22}^{-1} A_{21}, \;\;\;\bar N_k:=N_{k, 11}-N_{k, 12} A_{22}^{-1} A_{21},
      \end{align*}     
where $x_r(0)=0$ and the time dependence in (\ref{romstateeq}) is omitted to shorten the notation. This straight forward ansatz is based on observations 
from the deterministic case ($N_k=H_i=0$), where $x_2$ represents the fast variables, i.e., $\dot x_2(t) \approx 0$ after a short time, see \cite{spa}.

This ansatz for stochastic systems might, however, be false, no matter how small the HSVs corresponding to $x_2$ are. Despite the fact that for the motivation, a maybe less convincing argument is used, 
this leads to a viable MOR method for which an error bound can be proved. An averaging principle would be a mathematically well-founded alternative to this naive approach. Averaging principles for
stochastic systems have for example been investigated in \cite{avp2,avp1}. A further strategy to derive a ROM in this context can be found in \cite{berglundgentz}.\smallskip

Moreover, notice that system (\ref{romstochstatebt}) is not a bilinear system anymore due to the quadratic term in the control $u$. This is an essential difference to the ROM proposed in \cite{hartmann}.

\subsection{Main result}

The work in this paper on SPA for system (\ref{controlsystemoriginal}) can be interpreted as a generalization of the deterministic bilinear case \cite{hartmann}. This extension builds a bridge between stochastic
linear systems and stochastic nonlinear systems such that SPA can possibly be applied to many more stochastic equations and applications.\smallskip

In this paper, we provide an alternative to \cite{redstochbil}, where BT was studied. We extend the work of \cite{hartmann} combined with a modification of the ROM and the choice of a new Gramian defined through 
(\ref{newgram2}). Based on this, we obtain an error bound that was not even available for the deterministic bilinear case. This is the main result of this paper and is formulated in the following theorem.
Its proof requires new techniques that cannot be found in the literature so far.
\begin{theorem}\label{mainthmintro}
Let $y$ be the output of the full model (\ref{controlsystemoriginal}) with $x(0)=0$ and $y_r$ be the output of the ROM (\ref{romstochstatebt}) with zero initial state. Then, for all $T>0$, 
it holds that \begin{align*}
  \left(\mathbb E\left\|y-y_r\right\|_{L^2_{T}}^2\right)^{\frac{1}{2}}\leq 2 (\tilde \sigma_{1}+\tilde \sigma_{2}+\ldots + \tilde \sigma_\nu)  \left\|u\right\|_{L^2_T}\exp\left(0.5 \left\|u^0\right\|_{L^2_T}^2\right),            
  \end{align*}
where $\tilde\sigma_{1}, \tilde\sigma_{2}, \ldots,\tilde\sigma_\nu$ are the distinct diagonal entries of $\Sigma_2=\diag(\sigma_{r+1},\ldots,\sigma_n)=\diag(\tilde\sigma_{1} I, \tilde\sigma_{2} I, \ldots, \tilde\sigma_\nu I)$
and $u^0=(u^0_1, \dots, u_m^0)^T$ is the control vector with components defined by $u_k^0 \equiv \begin{cases}   0  & \text{if }N_k = 0,\\   u_k & \text{else}. \end{cases}$
\end{theorem}\smallskip

Theorem \ref{mainthmintro} is proved in Section \ref{proofmainthm}. We observe that an exponential term enters the bound in Theorem \ref{mainthmintro} which is due to the bilinearity in the drift. Setting 
$N_k=0$ for all $k=1, \ldots, m$ the exponential becomes a one which is the bound of the stochastic linear case \cite{redmannspa2}. The result in Theorem \ref{mainthmintro} tells us that 
the ROM (\ref{romstochstatebt}) yields a very good approximation if the truncated HSVs (diagonal entries of $\Sigma_2$) are small and the vector $u^0$ of control components with a non-zero 
$N_k$ is not too large. The exponential in the error bound can be an indicator that SPA performs badly if $u^0$ is very large. \smallskip

The remainder of the paper deals with the proof of Theorem \ref{mainthmintro}.

\section{$L^2$-error bound for SPA}\label{errorboundsBT}

The proof of the error bound in Theorem \ref{mainthmintro} is divided into two parts. We first investigate the error that we encounter by removing the smallest HSV from the system in Section \ref{sectionremovesmallhsv}. 
In this reduction step, the structure from the full model (\ref{controlsystemoriginal}) to the ROM (\ref{romstochstatebt}) changes. Therefore, when removing the other HSVs from the system, another case needs to be studied 
in Section \ref{secneigeb}. There, an error bound between two ROM is achieved which are neighboring, i.e., the larger ROM has exactly one HSV more than the smaller one. The results of Sections 
\ref{sectionremovesmallhsv} and \ref{secneigeb} are then combined in Section \ref{proofmainthm} in order to prove the general error bound.\smallskip

For simplicity, let us from now on assume that system (\ref{controlsystemoriginal}) is already balanced and has a zero initial condition ($x_0=0$). Thus, (\ref{newgram2}) and (\ref{gengenlyapobs}) become
\begin{align}\label{balancedreach}
 A^T \Sigma^{-1}+\Sigma^{-1}A+\sum_{k=1}^m N_k^T \Sigma^{-1} N_k + \sum_{i, j=1}^v H_i^T \Sigma^{-1} H_j k_{i j} &\leq -\Sigma^{-1}BB^T \Sigma^{-1},\\ \label{balancedobserve}
 A^T \Sigma+\Sigma A+\sum_{k=1}^m N_k^T \Sigma N_k +\sum_{i, j=1}^v H_i^T \Sigma H_j k_{ij} &\leq -C^T C,
                                       \end{align}
i.e., $P=Q=\Sigma=\diag(\sigma_1, \ldots, \sigma_n)>0$. 

\subsection{Error bound of removing the smallest HSV}\label{sectionremovesmallhsv}

We introduce the variable  $x_{\mp} =\smat x_1-x_r \\ x_2+A_{22}^{-1}(A_{21} x_r + B_2u)\srix$ since the corresponding output 
\begin{align}\label{xminusoutput}
    y_{\mp}(t)&=Cx_{\mp}(t)=Cx(t)-\bar C x_r(t) - \bar D u = y(t)-y_r(t), \;\;\;t\geq 0,
            \end{align}
is the output error between the full and the reduced system (\ref{romstochstatebt}). We aim to find an equation for $x_{\mp}$. This is done through the state variable
$x_-=\smat x_1-x_r \\ x_2\srix$. The differential $d(x_1-x_r)$  is obtained by subtracting the state equation (\ref{romstateeq}) of the reduced system from the first $r$ rows of  (\ref{stateeq}). The corresponding 
right side is then rewritten using $x_{\mp}$. Moreover, the right side of the differential of $x_2$, compare with the last $n-r$ rows of (\ref{stateeq}), is also formulated with the help of $x_{\mp}$. This results in
\begin{align}\label{statexminus}
             dx_-&=[Ax_{\mp}+ \smat 0 \\ c_0\srix+\sum_{k=1}^m N_{k} x_{\mp} u_k]dt +\sum_{i=1}^v [H_{i} x_{\mp} + \smat 0 \\ c_i\srix] dM_i,
            \end{align}
where $c_0(t):=\sum_{k=1}^m [ N_{k, 21}x_r(t)  -N_{k, 22}A_{22}^{-1}(A_{21} x_r(t) + B_2 u(t))] u_k(t)$ and $c_i(t):= H_{i, 21} x_r(t)-H_{i, 22}A_{22}^{-1}(A_{21} x_r(t) + B_2 u(t))$ for $i=1, \ldots, v$.\smallskip

We furthermore introduce the reverse state to $x_{\mp}$ in terms of the signs. This is 
$x_{\pm}=\smat x_1+x_r \\ x_2-A_{22}^{-1}(A_{21} x_r + B_2u)\srix$. Using the state $x_+=\smat x_1+x_r \\ x_2\srix$, with a differential obtained by combining (\ref{stateeq}) and (\ref{romstateeq}) again, and expressing 
its right side with $x_{\pm}$, we have 
\begin{align}\label{xplus}
             dx_+=[Ax_{\pm}+2 B u- \smat 0 \\ c_0\srix +\sum_{k=1}^m N_{k} x_{\pm} u_k]dt +\sum_{i=1}^v [H_{i} x_{\pm} - \smat 0 \\ c_i\srix] dM_i.
            \end{align}
We will see that the proof of the error bound can be reduced to the task of finding suitable estimates for $\mathbb E[x_-^T(t) \Sigma x_-(t)]$ and $\mathbb E[x_+^T(t) \Sigma^{-1} x_+(t)]$. This idea was 
also used to determine an error bound for BT \cite{redstochbil}. However, the proof for SPA requires different techniques to find the estimates.
  \begin{theorem}\label{mainthm}
Let $y$ be the output of the full model (\ref{controlsystemoriginal}) with $x(0)=0$, $y_r$ be the output of the ROM (\ref{romstochstatebt}) with $x_{r}(0)=0$ and $\Sigma_2=\sigma I$, $\sigma>0$, in (\ref{partitionfullmodel}).
Then, it holds that \begin{align*}
  \left(\mathbb E\left\|y-y_r\right\|_{L^2_{T}}^2\right)^{\frac{1}{2}}\leq 2 \sigma  \left\|u\right\|_{L^2_T}\exp\left(0.5 \left\|u^0\right\|_{L^2_T}^2\right).            
  \end{align*}
\begin{proof}
We derive a suitable upper bound for $\mathbb E[x_-^T(t) \Sigma x_-(t)]$ first applying Ito's formula. Hence, Lemma \ref{lemstochdiff} and Equation (\ref{statexminus}) yield
\begin{align}\label{productruleapplied}
\mathbb E\left[x_-^T(t)\Sigma x_-(t)\right]=&2 \int_0^t\mathbb E\left[x_-^T\Sigma\left(Ax_{\mp}+\sum_{k=1}^m (N_{k} x_{\mp} u_k) + \smat 0 \\ c_0\srix \right)\right]ds\\ \nonumber
&+ \int_0^t\sum_{i, j=1}^v \mathbb E\left[\left(H_{i} x_{\mp} + \smat 0 \\ c_i\srix\right)^T\Sigma\left(H_{j} x_{\mp} + \smat 0 \\ c_j\srix\right)\right]k_{ij} ds.
\end{align}
We find an estimate for the terms related to $N_k$, that is
\begin{align}\label{Nkestimate}
\sum_{k=1}^m 2 x_-^T(s)\Sigma N_{k} x_{\mp}(s) u_k(s)&=\sum_{k=1}^m 2\left\langle  \Sigma^{\frac{1}{2}} x_-(s)u_k(s), \Sigma^{\frac{1}{2}} N_k x_{\mp}(s) \right\rangle_2 \\ \nonumber
&\leq   \sum_{k=1}^m \left\|  \Sigma^{\frac{1}{2}} x_-(s)u^0_k(s)\right\|_2^2  + \left\|\Sigma^{\frac{1}{2}} N_k x_{\mp}(s)\right\|_2^2\\ \nonumber
&=x_-^T(s) \Sigma x_-(s)  \left\|u^0(s)\right\|_{2}^2 +\sum_{k=1}^m  x_{\mp}^T(s) N_k^T\Sigma N_{k} x_{\mp}(s),
\end{align}
where $u^0$ is defined as in Theorem \ref{mainthmintro}. Moreover, adding a zero, we rewrite \begin{align}\label{Aestimate}
2 x_-^T(s) \Sigma Ax_{\mp}(s) &=2 x_{\mp}^T(s) \Sigma Ax_{\mp}(s)- 2 \smat 0 \\ h(s)\srix^T \Sigma Ax_{\mp}(s) \\ \nonumber
& = x_{\mp}^T(s) (A^T \Sigma + \Sigma A) x_{\mp}(s) - 2 \smat 0 \\ h(s)\srix^T \Sigma Ax_{\mp}(s),                                   
                                    \end{align}
where $h(s)=A_{22}^{-1}(A_{21} x_r(s) + B_2u(s))$
With (\ref{Nkestimate}) and (\ref{Aestimate}), (\ref{productruleapplied}) becomes 
{\allowdisplaybreaks\begin{align}\nonumber
\mathbb E\left[x_-^T(t)\Sigma x_-(t)\right]\leq& \mathbb E\int_0^t x_{\mp}^T\left(A^T \Sigma+\Sigma A+\sum_{k=1}^m  N_k^T\Sigma N_{k}+\sum_{i, j=1}^v H^T_{i}\Sigma H_{j}k_{ij}\right)x_{\mp} ds\\ \label{insetobserveequation}
& +\mathbb E\int_0^t 2 x_-^T\Sigma \smat 0 \\ c_0\srix + \sum_{i, j=1}^v \left(2 H_{i} x_{\mp} + \smat 0 \\ c_i\srix\right)^T\Sigma\smat 0 \\ c_j\srix k_{ij}  ds\\
&+\int_0^t \mathbb E\left[x_-^T \Sigma x_-\right]  \left\|u^0\right\|_{2}^2 ds - \mathbb E \int_0^t 2 \smat 0 \\ h\srix^T \Sigma Ax_{\mp}ds.\nonumber
\end{align}}
Taking the partitions of $x_-$ and $\Sigma$ into account, we see that $x_-^T\Sigma \smat 0 \\ c_0\srix=x_2^T \Sigma_2 c_0$. Furthermore, the partitions of $x_{\mp}$ and $H_i$ yield 
\begin{align}\label{relforH}
&\left(2 H_{i} x_{\mp} + \smat 0 \\ c_i\srix\right)^T\Sigma\smat 0 \\ c_j\srix =\left(2 H_{i} x_{\mp} + \smat 0 \\ c_i\srix\right)^T\smat 0 \\ \Sigma_2 c_j\srix\\ \nonumber
&=\left(2 \smat H_{i, 21} & H_{i, 22}\srix (x - \smat x_r \\ -h\srix) +  c_i\right)^T \Sigma_2 c_j = \left(2 \smat H_{i, 21} & H_{i, 22}\srix x -  c_i\right)^T \Sigma_2 c_j,
\end{align}
since $\smat H_{i, 21} & H_{i, 22}\srix  \smat x_r \\ -h\srix=c_i$. Using the partition of $A$, it holds that \begin{align}\label{relforA}
-2 \smat 0 \\ h\srix^T \Sigma Ax_{\mp}&=-2 \smat 0 & h^T\Sigma_2 \srix Ax_{\mp}=-2 h^T\Sigma_2  \smat A_{21} & A_{22}\srix (x+\smat -x_r \\ h\srix)\\ \nonumber
&=-2 h^T\Sigma_2  (\smat A_{21} & A_{22}\srix x + B_2u),
\end{align}
because $\smat A_{21} & A_{22}\srix \smat -x_r \\ h\srix=B_2u$. We insert (\ref{balancedobserve}) and (\ref{xminusoutput}) into inequality (\ref{insetobserveequation}) and exploit the
relations in (\ref{relforH}) and (\ref{relforA}). Hence, 
\begin{align*}
\mathbb E\left[x_-^T(t)\Sigma x_-(t)\right]\leq& - \mathbb E\left\|y-y_r\right\|^2_{L^2_{t}}+\int_0^t \mathbb E\left[x_-^T \Sigma x_-\right]  \left\|u^0\right\|_{2}^2 ds\\ 
& +\mathbb E\int_0^t 2 x_2^T \Sigma_2 c_0 + \sum_{i, j=1}^v \left(2 \smat H_{i, 21} & H_{i, 22}\srix x -  c_i\right)^T \Sigma_2 c_j k_{ij}  ds\\
& - \mathbb E \int_0^t 2 h^T\Sigma_2  (\smat A_{21} & A_{22}\srix x + B_2u) ds.
\end{align*}
We define the function $\alpha_-(t):=\mathbb E\int_0^t 2 x_2^T \Sigma_2 c_0 + \sum_{i, j=1}^v \left(2 \smat H_{i, 21} & H_{i, 22}\srix x -  c_i\right)^T \Sigma_2 c_j k_{ij}  ds
- \mathbb E \int_0^t 2 h^T\Sigma_2  (\smat A_{21} & A_{22}\srix x + B_2u) ds$ and apply Lemma \ref{gronwall} implying
\begin{align*}
\mathbb E\left[x_-^T(t)\Sigma x_-(t)\right]\leq& \alpha_-(t)- \mathbb E\left\|y-y_r\right\|_{L^2_{t}}^2\\
& +\int_0^t (\alpha_-(s) - \mathbb E\left\|y-y_r\right\|_{L^2_{s}}^2) \left\|u^0(s)\right\|_{2}^2 \exp\left(\int_s^t \left\|u^0(w)\right\|_{2}^2 dw\right) ds.
\end{align*}
Since $\Sigma$ is positive definite, we obtain and upper bound for the output error by \begin{align*}
\mathbb E\left\|y-y_r\right\|_{L^2_{t}}^2\leq \alpha_-(t) +\int_0^t \alpha_-(s) \left\|u^0(s)\right\|_{2}^2 \exp\left(\int_s^t \left\|u^0(w)\right\|_{2}^2 dw\right) ds.
\end{align*}
Defining the term $\alpha_+(t):=\mathbb E\int_0^t 2 x_2^T \Sigma_2^{-1} c_0 + \sum_{i, j=1}^v \left(2 \smat H_{i, 21} & H_{i, 22}\srix x -  c_i\right)^T \Sigma_2^{-1} c_j k_{ij} ds
- \mathbb E \int_0^t 2 h^T\Sigma_2^{-1}  (\smat A_{21} & A_{22}\srix x + B_2u) ds$
and exploiting the assumption that $\Sigma_2=\sigma I$, leads to
\begin{align}\label{firstbound}
\mathbb E\left\|y-y_r\right\|_{L^2_{t}}^2\leq \sigma^2\left[\alpha_+(t) +\int_0^t \alpha_+(s) \left\|u^0(s)\right\|_{2}^2 \exp\left(\int_s^t \left\|u^0(w)\right\|_{2}^2 dw\right) ds\right].
\end{align}
The remaining step is to find a bound for the right side of (\ref{firstbound}) that does not depend on $\alpha_+$ anymore. For that reason, 
a bound for the expression $\mathbb E[x_+^T(t) \Sigma^{-1} x_+(t)]$ is derived next using Ito's lemma again. From (\ref{xplus}) and Lemma \ref{lemstochdiff}, we obtain
\begin{align}\label{productruleappliedplus}
\mathbb E\left[x_+^T(t)\Sigma^{-1} x_+(t)\right]=&2 \int_0^t\mathbb E\left[x_+^T\Sigma^{-1}\left(Ax_{\pm}+2 Bu+\sum_{k=1}^m (N_{k} x_{\pm} u_k) - \smat 0 \\ c_0\srix \right)\right]ds\\ \nonumber
&+ \int_0^t\sum_{i, j=1}^v \mathbb E\left[\left(H_{i} x_{\pm} - \smat 0 \\ c_i\srix\right)^T\Sigma^{-1}\left(H_{j} x_{\pm} - \smat 0 \\ c_j\srix\right)\right]k_{ij} ds.
\end{align}
Analogously to (\ref{Nkestimate}), it holds that \begin{align*}
\sum_{k=1}^m 2 x_+^T(s)\Sigma^{-1} N_{k} x_{\pm}(s) u_k(s) \leq  x_+^T(s) \Sigma^{-1} x_+(s)  \left\|u^0(s)\right\|_{2}^2 +\sum_{k=1}^m  x_{\pm}^T(s) N_k^T\Sigma^{-1} N_{k} x_{\pm}(s).
\end{align*}
Additionally, we rearrange the term related to $A$ as follows \begin{align*}
  2x_+^T(s) \Sigma^{-1} Ax_{\pm}(s) &= 2x_{\pm}^T(s) \Sigma^{-1} Ax_{\pm}(s)+ 2\smat 0 \\ h(s)\srix^T \Sigma^{-1} Ax_{\pm}(s)  \\
  &=  x_{\pm}^T(s)(A^T \Sigma^{-1} + \Sigma^{-1} A)x_{\pm}(s)+ 2\smat 0 \\ h(s)\srix^T \Sigma^{-1} Ax_{\pm}(s).                                                       
                                                           \end{align*}
Moreover, we have \begin{align*}
 4 x_+^T(s)\Sigma^{-1} Bu(s) = 4 x_{\pm}^T(s)\Sigma^{-1} Bu(s)+ 4\smat 0 \\ h(s)\srix^T \Sigma^{-1} Bu(s).                                
                                                           \end{align*}
We plug in the above results into (\ref{productruleappliedplus}) which gives us 
\begin{align}\nonumber
&\mathbb E\left[x_+^T(t)\Sigma^{-1} x_+(t)\right]\\ \nonumber&\leq \mathbb E\int_0^t x_{\pm}^T\left(A^T \Sigma^{-1}+\Sigma^{-1} A+\sum_{k=1}^m  N_k^T\Sigma^{-1} N_{k}+\sum_{i, j=1}^v H^T_{i}\Sigma^{-1} H_{j}k_{ij}\right)x_{\pm} ds\\ \label{insetobserveequation2}
&\quad -\mathbb E\int_0^t 2 x_+^T\Sigma^{-1} \smat 0 \\ c_0\srix + \sum_{i, j=1}^v \left(2 H_{i} x_{\pm} - \smat 0 \\ c_i\srix\right)^T\Sigma^{-1}\smat 0 \\ c_j\srix k_{ij}  ds\\
&\quad+ \mathbb E \int_0^t 2\smat 0 \\ h\srix^T \Sigma^{-1} (Ax_{\pm}+2 Bu) ds + \mathbb E \int_0^t 4 x_{\pm}^T\Sigma^{-1} Bu ds\nonumber \\
&\quad+\int_0^t \mathbb E\left[x_+^T \Sigma^{-1} x_+\right]  \left\|u^0\right\|_{2}^2 ds. \nonumber
\end{align}
From inequality (\ref{balancedreach}) and the Schur complement condition on definiteness, it follows that\begin{align}\label{schurposdef}
 \mat{cc}\hspace{-0.15cm} A^T \Sigma^{-1}\hspace{-0.025cm}+\Sigma^{-1}\hspace{-0.025cm}A+\hspace{-0.025cm}\sum_{k=1}^m N_k^T \Sigma^{-1} N_k + \hspace{-0.025cm}\sum_{i, j=1}^v H_i^T \Sigma^{-1} H_j k_{i j}  & \Sigma^{-1}B\\
 B^T \Sigma^{-1}& -I\rix\leq 0.
                                       \end{align}
We multiply (\ref{schurposdef}) with $\smat x_{\pm} \\ 2u\srix^T$ from the left and with $\smat x_{\pm}\\ 2u\srix$ from the right. Hence, 
\begin{align}\label{ugrosserlyap}
&4 \left\|u\right\|_{2}^2\geq \\ \nonumber
& x_{\pm}^T\left(A^T \Sigma^{-1}+\Sigma^{-1} A+\sum_{k=1}^m  N_k^T\Sigma^{-1} N_{k}+\sum_{i, j=1}^v H^T_{i}\Sigma^{-1} H_{j}k_{ij}\right)x_{\pm}+4x_{\pm}^T\Sigma^{-1} Bu.
\end{align}
Applying this result to (\ref{insetobserveequation2}) yields 
{\allowdisplaybreaks\begin{align}\label{insetobserveequationbla}
\mathbb E\left[x_+^T(t)\Sigma^{-1} x_+(t)\right] \leq & 4 \left\|u\right\|_{L^2_t}^2+\int_0^t \mathbb E\left[x_+^T \Sigma^{-1} x_+\right]  \left\|u^0\right\|_{2}^2 ds\\\nonumber
& +\mathbb E \int_0^t 2\smat 0 \\ h\srix^T \Sigma^{-1} (Ax_{\pm}+2 Bu) ds\\ \nonumber
& -\mathbb E\int_0^t 2 x_+^T\Sigma^{-1} \smat 0 \\ c_0\srix + \sum_{i, j=1}^v \left(2 H_{i} x_{\pm} - \smat 0 \\ c_i\srix\right)^T\Sigma^{-1}\smat 0 \\ c_j\srix k_{ij}  ds.
\end{align}}
We first of all see that $x_+^T\Sigma^{-1} \smat 0 \\ c_0\srix=x_2^T\Sigma_2^{-1}c_0$ using the partitions of $x_+$ and $\Sigma$.  
 With the partition of $H_i$, we moreover have 
\begin{align*}
&\left(2 H_{i} x_{\pm} - \smat 0 \\ c_i\srix\right)^T\Sigma^{-1}\smat 0 \\ c_j\srix =\left(2 H_{i} x_{\pm} - \smat 0 \\ c_i\srix\right)^T\smat 0 \\ \Sigma_2^{-1} c_j\srix\\
&=\left(2 \smat H_{i, 21} & H_{i, 22}\srix (x + \smat x_r \\ -h\srix) -  c_i\right)^T \Sigma_2^{-1} c_j = \left(2 \smat H_{i, 21} & H_{i, 22}\srix x +  c_i\right)^T \Sigma_2^{-1} c_j.
\end{align*}
In addition, it holds that \begin{align*}
&2\smat 0 \\ h\srix^T \Sigma^{-1} (Ax_{\pm}+2 Bu) = 2\smat 0 & h^T\Sigma_2^{-1}\srix  (Ax_{\pm}+2 Bu)\\
&= 2 h^T\Sigma_2^{-1} (\smat A_{21} & A_{22}\srix (x+\smat x_r \\ -h\srix)+2 B_2u)   =   2 h^T\Sigma_2^{-1} (\smat A_{21} & A_{22}\srix x+ B_2u)                        
                           \end{align*}
Plugging the above relations into (\ref{insetobserveequationbla}) leads to
\begin{align}\label{insetobserveequationblabla}
\mathbb E\left[x_+^T(t)\Sigma^{-1} x_+(t)\right] \leq & 4\left\|u\right\|_{L^2_t}^2+\int_0^t \mathbb E\left[x_+^T \Sigma^{-1} x_+\right]  \left\|u^0\right\|_{2}^2 ds\\\nonumber
&+\mathbb E\int_0^t 2 h^T\Sigma_2^{-1} (\smat A_{21} & A_{22}\srix x+ B_2u) ds\\ \nonumber
& -\mathbb E\int_0^t 2 x_2^T\Sigma_2^{-1} c_0 + \sum_{i, j=1}^v \left(2 \smat H_{i, 21} & H_{i, 22}\srix x + c_i\right)^T\Sigma_2^{-1} c_j k_{ij} ds.
\end{align}
We add $2\mathbb E\int_0^t \sum_{i, j=1}^v c_i^T\Sigma_2^{-1} c_j k_{ij} ds$ to the right side of (\ref{insetobserveequationblabla}) and preserve the inequality since this term is a nonnegative due to Lemma \ref{proppossemidef}. 
This results in
\begin{align*}
\mathbb E\left[x_+^T(t)\Sigma^{-1} x_+(t)\right] \leq 4 \left\|u\right\|_{L^2_t}^2-\alpha_+(t)+\int_0^t \mathbb E\left[x_+^T(s) \Sigma^{-1} x_+(s)\right]  \left\|u^0(s)\right\|_{2}^2 ds.
\end{align*}
Gronwall's inequality in Lemma \ref{gronwall} yields \begin{align}\label{keineahnung}
&\mathbb E\left[x_+^T(t)\Sigma^{-1} x_+(t)\right] \\ \nonumber &\leq  4\left\|u\right\|_{L^2_t}^2-\alpha_+(t) 
+\int_0^t (4 \left\|u\right\|_{L^2_s}^2-\alpha_+(s))  \left\|u^0(s)\right\|_{2}^2 \exp\left(\int_s^t \left\|u^0(w)\right\|_{2}^2 dw\right) ds.
\end{align}
We find an estimate for the following expression: \begin{align}\label{keineahnung2}
&\int_0^t \left\|u\right\|_{L^2_s}^2 \left\|u^0(s)\right\|_2^2 \exp\left(\int_s^t \left\|u^0(w)\right\|_2^2dw\right) ds\\ \nonumber
&\leq \left\|u\right\|_{L^2_t}^2 \left[-\exp\left(\int_s^t \left\|u^0(w)\right\|_2^2dw\right)\right]_{s=0}^t\\
&=\left\|u\right\|_{L^2_t}^2\left(\exp\left(\int_0^t \left\|u^0(s)\right\|_2^2ds\right)-1\right). \nonumber
\end{align}
Combining (\ref{keineahnung}) with (\ref{keineahnung2}), we obtain
\begin{align}\label{analoggronwallest}
&\alpha_+(t)+\int_0^t \alpha_+(s)  \left\|u^0(s)\right\|_{2}^2 \exp\left(\int_s^t \left\|u^0(w)\right\|_{2}^2 dw\right) ds\\ \nonumber
&\leq 4\left\|u\right\|_{L^2_t}^2 \exp\left(\int_0^t \left\|u^0(s)\right\|_2^2ds\right).
\end{align}
Comparing this result with (\ref{firstbound}) implies 
 \begin{align}\label{secondbound}
\left(\mathbb E\left\|y-y_r\right\|_{L^2_{t}}^2\right)^{\frac{1}{2}}\leq 2\sigma \left\|u\right\|_{L^2_t} \exp\left(0.5 \left\|u^0\right\|_{L^2_t}^2\right).
\end{align}
\end{proof}
\end{theorem}\\
We proceed with the study of an error bound between two ROM that are neighboring.
\subsection{Error bound for neighboring ROMs}\label{secneigeb}
In this section, we investigate the output error between two ROMs, in which the larger ROM has exactly one HSV than the smaller one. This concept of neighboring ROMs was first introduced in \cite{redmannspa2} but in the much simpler 
stochastic linear setting.\smallskip

The reader might wonder why a second case is considered besides the one in Section 
\ref{sectionremovesmallhsv} since one might just start with a full model that has the same structure as the ROM (\ref{romstochstatebt}). The reason is that is not clear how the Gramians need to be chosen for (\ref{romstochstatebt}). 
In order to investigate the error between two ROMs by SPA, a finer partition than the one in (\ref{partvonorimodel}) is required. We partition the matrices of the balanced full system (\ref{controlsystemoriginal}) as follows:
\begin{subequations}\label{finerpartdef}
\begin{align}
  A=\smat{A}_{11}&{A}_{12}&A_{13}\\ 
{A}_{21}&{A}_{22}&A_{23}\\
{A}_{31}&{A}_{32}&A_{33}\srix,\quad B=\smat B_1 \\ B_2\\B_3\srix,\quad C= \smat C_1 & C_2& C_3\srix,\\
H_i=\smat{H}_{i, 11}&{H}_{i, 12}&{H}_{i, 13}\\ 
{H}_{i, 21}&{H}_{i, 22}&{H}_{i, 23}\\
{H}_{i, 31}&{H}_{i, 32}&{H}_{i, 33}\srix,\quad N_k=\smat{N}_{k, 11}&{N}_{k, 12}&{N}_{k, 13}\\ 
{N}_{k, 21}&{N}_{k, 22}&{N}_{k, 23}\\
{N}_{k, 31}&{N}_{k, 32}&{N}_{k, 33}\srix.
            \end{align}
            \end{subequations}
The partitioned balanced solution to (\ref{stateeq}) and the Gramians are then of the form \begin{align}\label{finepartsig}
 x=\smat x_{1}\\ 
 x_{2}\\
 x_{3}\srix\;  \text{and} \;       
\Sigma=\smat \Sigma_{1}& & \\ 
 &\Sigma_{2}& \\
 & &\Sigma_{3}\srix.                                                   
                                                     \end{align}
We introduce the ROM of truncating $\Sigma_3$ first. According to the procedure described in Section \ref{settingstochstabgen}, the reduced system is obtained by setting $dx_3$ equal to zero, neglecting 
the bilinear and the diffusion term in this equation. The solution $\tilde x_3$ of the resulting algebraic constraint is an approximation for $x_3$. 
One can solve for this approximating variable and obtains $\tilde x_3=-A_{33}^{-1}(A_{31}x_1+A_{32}x_2+B_3u)$. Inserting this result for $x_3$ in the equations for $x_1$, $x_2$ and into the output equation (\ref{originalobserveq})
leads to 
\begin{subequations}\label{romnegsig3}
 \begin{align}\label{stateromsig3}
d\smat x_1 \\ x_2\srix&=\left[\hat A\smat x_1 \\ x_2\\\tilde x_3\srix+\hat B u+\sum_{k=1}^m \hat N_k \smat x_1 \\ 
x_2\\  \tilde x_3\srix u_k\right]dt +\sum_{i=1}^v \hat H_i \smat x_1 \\ x_2\\ \tilde x_3\srix dM_i, \\ \label{outromsig3}
             \bar y(t)&= C \smat x_1(t) \\ x_2(t)\\ \tilde x_3(t)\srix, \;\;\;t\geq 0,
            \end{align}
            \end{subequations}
where $\smat x_1(0) \\ x_2(0)\srix=\smat 0 \\ 0\srix$ and 
\begin{align*}
\hat A=\smat{A}_{11}&{A}_{12}&A_{13}\\ 
{A}_{21}&{A}_{22}&A_{23}\srix,\; \hat B=\smat B_1 \\ B_2\srix, \; \hat H_i=\smat{H}_{i, 11}&{H}_{i, 12}&{H}_{i, 13}\\ 
{H}_{i, 21}&{H}_{i, 22}&{H}_{i, 23}\srix,\; \hat N_k=\smat{N}_{k, 11}&{N}_{k, 12}&{N}_{k, 13}\\ 
{N}_{k, 21}&{N}_{k, 22}&{N}_{k, 23}\srix.
            \end{align*}
We aim to determine the error between this ROM and the reduced system of neglecting $\Sigma_2$  and $\Sigma_3$. This is 
\begin{subequations}\label{romnegsig2sig3}
 \begin{align}\label{stateromsig2sig3}
dx_r&=\left[\hat A_r \smat x_r \\ -h_1\\ -h_2\srix+B_1 u+\sum_{k=1}^m \hat N_{r, k} \smat x_r \\ 
-h_1\\ -h_2\srix u_k\right]dt + \sum_{i=1}^v \hat H_{r, i} \smat x_r \\ -h_1\\ -h_2\srix dM_i, \\ \label{outromsig2sig3}
             \bar y_r(t)&=\smat C_1 & C_2 &C_3 \srix\smat x_r(t) \\ -h_1(t)\\ -h_2(t)\srix, \;\;\;t\geq 0,
            \end{align}
\end{subequations}
where $x_r(0)=0$,
\begin{align*}
\hat A_r=\smat{A}_{11}&{A}_{12}&A_{13}\srix,\; \hat H_{r, i}=\smat{H}_{i, 11}&{H}_{i, 12}&{H}_{i, 13}\srix,\; \hat N_{r, k}=\smat{N}_{k, 11}&{N}_{k, 12}&{N}_{k, 13}\srix   \end{align*}
and we define
\begin{align}\label{inverserep}
   h(t)=\smat h_1(t) \\ h_2(t)\srix= \smat A_{22}& A_{23}\\ A_{32}& A_{33}\srix^{-1}   
   \left(\smat {A}_{21}\\{A}_{31}\srix x_r(t)+\smat {B}_{2}\\{B}_{3}\srix u(t)\right).
                        \end{align}
In order to find a bound for the error between (\ref{outromsig3}) and (\ref{outromsig2sig3}), state variables analogously to $x_\mp$ and $x_\pm$ in Section \ref{sectionremovesmallhsv} are constructed in the following and 
corresponding equations are derived. For simplicity, we use a similar notation again and define
\begin{align*}
\hat x_\mp=\smat x_1-x_r \\ x_2+h_1\\ \tilde x_3+h_2\srix\; \text{and}\; \hat x_\pm=\smat x_1+x_r \\ x_2-h_1\\ \tilde x_3-h_2\srix.
\end{align*}
One can see that these states are obtained by combining the states appearing on the right sides of (\ref{stateromsig3}) and (\ref{stateromsig2sig3}). Furthermore, the output of $\hat x_\mp$ leads to the output error
\begin{align}\label{outputerrorbrom}
C\hat x_\mp(t)=  \bar y(t) -  \bar y_r(t), \quad t\geq 0,
\end{align}
which is a direct consequence of (\ref{outromsig3}) and (\ref{outromsig2sig3}).\smallskip

Now, we find the differential equations for $\hat x_\mp$ and $\hat x_\mp$. Using (\ref{inverserep}), we find that
\begin{align}\nonumber
 \smat {A}_{21}&{A}_{22}& A_{23}\\{A}_{31}&{A}_{32}& A_{33}  \srix \smat x_r \\ 
-h_1\\-h_2\srix &= \smat {A}_{21}\\ {A}_{31}\srix x_r -\smat {A}_{22}& A_{23}\\{A}_{32}& A_{33}  \srix h\\
&= \smat {A}_{21}\\ {A}_{31}\srix x_r -\smat {A}_{22}& A_{23}\\{A}_{32}& A_{33}  \srix \smat {A}_{22}& A_{23}\\{A}_{32}& A_{33}  \srix^{-1} \left(\smat {A}_{21}\\{A}_{31}\srix x_r+\smat {B}_{2}\\{B}_{3}\srix u\right)\nonumber \\
& = -\smat {B}_{2}\\{B}_{3}\srix u.\label{reltobeusef}
                         \end{align}            
Applying the first line of (\ref{reltobeusef}), we obtain the following equation \begin{align}\nonumber
d  0 &=\left[\smat {A}_{21}&{A}_{22}& A_{23}\srix\smat x_r \\ -h_1\\ -h_2\srix+ {B}_{2} u -  \hat c_0+\sum_{k=1}^m \smat{N}_{k, 21}&{N}_{k, 22}&{N}_{k, 23} \srix \smat x_r \\ 
-h_1\\ -h_2\srix u_k\right]dt\\ \label{zerosde}
&\quad+\sum_{i=1}^v \left[\smat{H}_{i, 21}&{H}_{i, 22}&{H}_{i, 23} \srix \smat x_r \\ -h_1\\ -h_2\srix - \hat c_i\right] dM_i
\end{align}
where $\hat c_0=\sum_{k=1}^m \smat{N}_{k, 21}&{N}_{k, 22}&{N}_{k, 23} \srix \smat x_r \\ -h_1\\ -h_2\srix u_k$ and 
$\hat c_i=\smat{H}_{i, 21}&{H}_{i, 22}&{H}_{i, 23} \srix \smat x_r \\ -h_1\\ -h_2\srix$
for $i=1, \ldots, v$. We supplement (\ref{stateromsig2sig3}) with (\ref{zerosde}) and combine this with (\ref{stateromsig3}). Hence, we obtain
\begin{align}\label{xminplu2}
d \hat x_-&=\left[\hat A \hat x_\mp+ \smat 0\\ \hat c_0\srix +\sum_{k=1}^m \hat {N}_{k} \hat x_\mp u_k\right]dt +\sum_{i=1}^v \left[\hat {H}_{i} \hat x_\mp + \smat 0\\ \hat c_i\srix\right] dM_i,
\end{align}
where $\hat x_- = \smat x_1 - x_r \\ x_2\srix$ and furthermore
\begin{align}\label{xplusmin2}
d \hat x_+&=\left[\hat A \hat x_\pm+ 2 \hat B u - \smat 0\\ \hat c_0\srix +\sum_{k=1}^m \hat {N}_{k} \hat x_\pm u_k\right]dt +\sum_{i=1}^v \left[\hat {H}_{i} \hat x_\pm - \smat 0\\ \hat c_i\srix\right] dM_i,
\end{align}
where $\hat x_+ = \smat x_1 + x_r \\ x_2\srix$. We now state the output error between the systems (\ref{romnegsig3}) and (\ref{romnegsig2sig3}) for the case that the ROM
are neighboring, i.e., the larger model has exactly one HSV more than the smaller one.
 \begin{theorem}\label{mainthm2}
Let $\bar y$ be the output of the ROM (\ref{romnegsig3}), $\bar y_r$ be the output of the ROM (\ref{romnegsig2sig3}) and $\Sigma_2=\sigma I$, $\sigma>0$, in (\ref{finepartsig}).
Then, it holds that \begin{align*}
  \left(\mathbb E\left\|\bar y-\bar y_r\right\|_{L^2_{T}}^2\right)^{\frac{1}{2}}\leq 2 \sigma  \left\|u\right\|_{L^2_T}\exp\left(0.5 \left\|u^0\right\|_{L^2_T}^2\right).            
  \end{align*}
\begin{proof}
We make use of equations (\ref{xminplu2}) and (\ref{xplusmin2}) in order to prove this bound. We set $\hat \Sigma= \smat \Sigma_1& \\ & \Sigma_2\srix$ as a submatrix of $\Sigma$ in (\ref{finepartsig}). 
Lemma \ref{lemstochdiff} now yields \begin{align}\label{productruleapplied2}
\mathbb E\left[\hat x_-^T(t)\hat \Sigma \hat x_-(t)\right]=&2 \int_0^t\mathbb E\left[\hat x_-^T\hat \Sigma\left(\hat A\hat x_{\mp}+\sum_{k=1}^m (\hat N_{k} \hat x_{\mp} u_k) + \smat 0 \\ \hat c_0\srix \right)\right]ds\\ \nonumber
&+ \int_0^t\sum_{i, j=1}^v \mathbb E\left[\left(\hat H_{i} \hat x_{\mp} + \smat 0 \\ \hat c_i\srix\right)^T\hat \Sigma\left(\hat H_{j} \hat x_{\mp} + \smat 0 \\ \hat c_j\srix\right)\right]k_{ij} ds.
\end{align}
We see that the right side of (\ref{productruleapplied2}) contains the submatrices $\hat A, \hat B, \hat H, \hat N$ and $\hat \Sigma$. In order to be able to refer to the full matrix inequality (\ref{balancedobserve}), we find 
upper bounds for certain terms in the following involving the full matrices $A, B, H, N$ and $\Sigma$. With the same estimate as in (\ref{Nkestimate}) and the control vector $u^0$ defined in Theorem \ref{mainthmintro}, we have
\begin{align*}
\sum_{k=1}^m 2 \hat x_-^T(s)\hat \Sigma \hat N_{k} \hat x_{\mp}(s) u_k(s)\leq   \hat x_-^T(s) \hat \Sigma \hat x_-(s)  \left\|u^0(s)\right\|_{2}^2 +\sum_{k=1}^m  \hat x_{\mp}^T(s) \hat N_k^T\hat \Sigma \hat N_{k} \hat x_{\mp}(s).
\end{align*}
Adding the term $\sum_{k=1}^m \left(\smat{N}_{k, 31}&{N}_{k, 32}&{N}_{k, 33} \srix \hat x_{\mp}(s)\right)^T \Sigma_3 \smat{N}_{k, 31}&{N}_{k, 32}&{N}_{k, 33} \srix \hat x_{\mp}(s)$ to the right side of this inequality results in
\begin{align}\label{Nkestimate2}
\sum_{k=1}^m 2 \hat x_-^T(s)\hat \Sigma \hat N_{k} \hat x_{\mp}(s) u_k(s)\leq   \hat x_-^T(s) \hat \Sigma \hat x_-(s)  \left\|u^0(s)\right\|_{2}^2 +\sum_{k=1}^m  \hat x_{\mp}^T(s) N_k^T \Sigma N_{k} \hat x_{\mp}(s).
\end{align}
Moreover, it holds that 
\begin{align*}
\hat x_{\mp}^T (A^T\Sigma+\Sigma A)\hat x_{\mp} &= 2 \hat x_{\mp}^T \Sigma A\hat x_{\mp} \\
&= 2 \smat x_1-x_r\\x_2+h_1 \srix^T \hat\Sigma \hat A\hat x_{\mp} + 2 (\tilde x_3+h_2)^T \Sigma_3\smat {A}_{31}&{A}_{32}& A_{33}\srix\hat x_{\mp}.
\end{align*}
We derive $\smat {A}_{31}&{A}_{32}& A_{33}\srix \smat x_{1}\\ x_{2}\\ \tilde x_{3}\srix=-B_3u$ by the definition of $\tilde x_3$. Moreover, it can be seen from the second line of (\ref{reltobeusef}) that
$\smat {A}_{31}&{A}_{32}& A_{33}\srix\hat x_{\mp}=0$. Hence, \begin{align}\label{estimateforea}
\hat x_{\mp}^T (A^T\Sigma+\Sigma A)\hat x_{\mp} = 2 \hat x_-^T \hat\Sigma \hat A\hat x_{\mp} + 2\smat 0\\ h_1\srix^T \hat \Sigma \hat A\hat x_{\mp}.
\end{align}
It remains to find a suitable upper bound related to the expression depending on $\hat H_i$. We first of all see that
\begin{align*}
 &\sum_{i, j=1}^v \left(\hat H_{i} \hat x_{\mp} + \smat 0 \\ \hat c_i\srix\right)^T\hat \Sigma\left(\hat H_{j} \hat x_{\mp} + \smat 0 \\ \hat c_j\srix\right)k_{ij}\\
 &= \hat x_{\mp}^T\sum_{i, j=1}^v \hat H^T_{i}\hat \Sigma \hat H_{j}k_{ij} \hat x_{\mp} + \sum_{i, j=1}^v \left(2 \hat H_{i} \hat x_{\mp} + \smat 0 \\ \hat c_i\srix\right)^T\hat \Sigma\smat 0 \\ \hat c_j\srix k_{ij}.
\end{align*}
The term $\sum_{i, j=1}^v \left(\smat{H}_{i, 31}&{H}_{i, 32}&{H}_{i, 33} \srix \hat x_{\mp}(s)\right)^T \Sigma_3 \smat{H}_{j, 31}&{H}_{j, 32}&{H}_{j, 33} \srix \hat x_{\mp}(s) k_{ij}$ 
is nonnegative through Lemma \ref{proppossemidef}. Adding this term to the right side of the above equation yields
\begin{align}\label{estimatereltoh}
 &\sum_{i, j=1}^v \left(\hat H_{i} \hat x_{\mp} + \smat 0 \\ \hat c_i\srix\right)^T\hat \Sigma\left(\hat H_{j} \hat x_{\mp} + \smat 0 \\ \hat c_j\srix\right)k_{ij}\\ \nonumber
 &\leq \hat x_{\mp}^T\sum_{i, j=1}^v H^T_{i} \Sigma H_{j}k_{ij} \hat x_{\mp} + \sum_{i, j=1}^v \left(2 \hat H_{i} \hat x_{\mp} + \smat 0 \\ \hat c_i\srix\right)^T\hat \Sigma\smat 0 \\ \hat c_j\srix k_{ij}.
\end{align}
Applying (\ref{Nkestimate2}), (\ref{estimateforea}) and (\ref{estimatereltoh}) to (\ref{productruleapplied2}), results in
{\allowdisplaybreaks\begin{align}\nonumber
\mathbb E\left[\hat x_-^T(t)\hat \Sigma \hat x_-(t)\right]\leq& \mathbb E\int_0^t \hat x_{\mp}^T\left(A^T \Sigma+\Sigma A+\sum_{k=1}^m  N_k^T\Sigma N_{k}+\sum_{i, j=1}^v H^T_{i}\Sigma H_{j}k_{ij}\right)\hat x_{\mp} ds\\
\label{insetobserveequation12}
& +\mathbb E\int_0^t 2 \hat x_-^T\hat \Sigma \smat 0 \\ \hat c_0\srix + \sum_{i, j=1}^v \left(2 \hat H_{i} \hat x_{\mp} + \smat 0 \\ \hat c_i\srix\right)^T \hat \Sigma\smat 0 \\ \hat c_j\srix k_{ij}  ds\\
&+\int_0^t \mathbb E\left[\hat x_-^T \hat \Sigma \hat x_-\right]  \left\|u^0\right\|_{2}^2 ds - \mathbb E \int_0^t 2 \smat 0 \\ h_1\srix^T \hat \Sigma \hat A \hat x_{\mp}ds.\nonumber
\end{align}}
Using that $\hat c_i=\smat{H}_{i, 21}&{H}_{i, 22}&{H}_{i, 23} \srix \smat x_r \\ -h_1\\ -h_2\srix$, we have
\begin{align}\nonumber
&\left(2 \hat H_{i} \hat x_{\mp} + \smat 0 \\ \hat c_i\srix\right)^T\hat \Sigma\smat 0 \\ \hat c_j\srix =\left(2 \hat H_{i} \hat x_{\mp} + \smat 0 \\ \hat c_i\srix\right)^T\smat 0 \\ \Sigma_2 \hat c_j\srix\\ \nonumber
&=\left(2 \smat H_{i, 21} & H_{i, 22} & H_{i, 23}\srix (\smat x_1\\ x_2 \\ \tilde x_3\srix - \smat x_r \\ -h_1\\-h_2\srix) +  \hat c_i\right)^T \Sigma_2 \hat c_j\\ \label{blaconfH}
 &= \left(2 \smat H_{i, 21} & H_{i, 22}& H_{i, 23}\srix  \smat x_1\\ x_2 \\ \tilde x_3\srix -  \hat c_i\right)^T \Sigma_2 \hat c_j.
\end{align}
It can be seen further that
\begin{align}\nonumber
-2\smat 0\\ h_1\srix^T \hat \Sigma \hat A\hat x_{\mp}&=-2 \smat 0 & h_1^T\Sigma_2 \srix \hat A\hat x_{\mp}=-2 h_1^T\Sigma_2  \smat A_{21} & A_{22}& A_{23}\srix (\smat x_1\\ x_2 \\ \tilde x_3\srix +\smat -x_r \\ h_1\\h_2\srix)\\
&=-2 h_1^T\Sigma_2  (\smat A_{21} & A_{22}& A_{23}\srix \smat x_1\\ x_2 \\ \tilde x_3\srix + B_2u)\label{relforA2}
\end{align}
taking the first line of (\ref{reltobeusef}) into account. Inserting (\ref{blaconfH}) and (\ref{relforA2}) into (\ref{insetobserveequation12}) and using the fact that 
$2 \hat x_-^T\hat \Sigma \smat 0 \\ \hat c_0\srix= 2 x_2 \Sigma_2 \hat  c_0  $ leads to 
{\allowdisplaybreaks\begin{align}\label{insetobserveequation122}
\mathbb E\left[\hat x_-^T(t)\hat \Sigma \hat x_-(t)\right]& \leq \int_0^t \mathbb E\left[\hat x_-^T \hat \Sigma \hat x_-\right]  \left\|u^0\right\|_{2}^2 ds +\hat \alpha_-(t)\\ \nonumber
&+\mathbb E\int_0^t \hat x_{\mp}^T\left(A^T \Sigma+\Sigma A+\sum_{k=1}^m  N_k^T\Sigma N_{k}+\sum_{i, j=1}^v H^T_{i}\Sigma H_{j}k_{ij}\right)\hat x_{\mp} ds,
\end{align}}
where we set $\hat \alpha_-(t):=\mathbb E\int_0^t 2 x_2^T \Sigma_2 \hat c_0 + \left(2 \smat H_{i, 21} & H_{i, 22}& H_{i, 23}\srix  \smat x_1\\ x_2 \\ \tilde x_3\srix -  \hat c_i\right)^T \Sigma_2 \hat c_j  ds
-\mathbb E\int_0^t 2 h_1^T\Sigma_2  (\smat A_{21} & A_{22}& A_{23}\srix \smat x_1\\ x_2 \\ \tilde x_3\srix + B_2u) ds$.
With (\ref{balancedobserve}) and (\ref{outputerrorbrom}), we obtain \begin{align*}
\mathbb E\left[\hat x_-^T(t)\hat \Sigma \hat x_-(t)\right]& \leq \int_0^t \mathbb E\left[\hat x_-^T \hat \Sigma \hat x_-\right]  \left\|u^0\right\|_{2}^2 ds +\hat \alpha_-(t)- \mathbb E\left\|\bar y-\bar y_r\right\|_{L^2_{t}}^2.
\end{align*}
Applying Lemma \ref{gronwall} to this inequality yields \begin{align*}
\mathbb E\left[\hat x_-^T(t)\hat \Sigma \hat x_-(t)\right]\leq& \hat \alpha_-(t)- \mathbb E\left\|\bar y-\bar y_r\right\|_{L^2_{t}}^2\\
& +\int_0^t \hat \alpha_-(s)\left\|u^0(s)\right\|_{2}^2 \exp\left(\int_s^t \left\|u^0(w)\right\|_{2}^2 dw\right) ds.
\end{align*}
Since the above left side of the inequality is positive, we obtain \begin{align*}
&\mathbb E\left\|\bar y - \bar y_r\right\|_{L^2_{t}}^2\\ &\leq \hat \alpha_-(t)
 +\int_0^t \hat \alpha_-(s)\left\|u^0(s)\right\|_{2}^2 \exp\left(\int_s^t \left\|u^0(w)\right\|_{2}^2 dw\right) ds.
\end{align*}
We exploit that $\Sigma_2=\sigma I$. Hence, we have \begin{align}\label{boundmitalpmin}
&\mathbb E\left\|\bar y - \bar y_r\right\|_{L^2_{t}}^2\\ \nonumber&\leq \sigma^2\left(\hat \alpha_+(t)
 +\int_0^t \hat \alpha_+(s)\left\|u^0(s)\right\|_{2}^2 \exp\left(\int_s^t \left\|u^0(w)\right\|_{2}^2 dw\right) ds\right),
\end{align}
where we set $\hat \alpha_+(t):=\mathbb E\int_0^t 2 x_2^T \Sigma_2^{-1} \hat c_0 + \left(2 \smat H_{i, 21} & H_{i, 22}& H_{i, 23}\srix  \smat x_1\\ x_2 \\ \tilde x_3\srix -  \hat c_i\right)^T \Sigma_2^{-1} \hat c_j  ds
-\mathbb E\int_0^t 2 h_1^T\Sigma_2^{-1}  (\smat A_{21} & A_{22}& A_{23}\srix \smat x_1\\ x_2 \\ \tilde x_3\srix + B_2u) ds$.
In order to find a suitable bound for the right side of (\ref{boundmitalpmin}),  Ito's lemma  is applied to $\mathbb E[\hat x_+^T(t) \hat \Sigma^{-1}\hat x_+(t)]$. Due to (\ref{xplusmin2}) and Lemma \ref{lemstochdiff}, we obtain
\begin{align}\label{productruleappliedplus22}
\mathbb E\left[\hat x_+^T(t)\hat \Sigma^{-1} \hat x_+(t)\right]=&2 \int_0^t\mathbb E\left[\hat x_+^T\hat \Sigma^{-1}\left(\hat A\hat x_{\pm}+2 \hat Bu+\sum_{k=1}^m (\hat N_{k} \hat x_{\pm} u_k) - \smat 0 \\ \hat c_0\srix \right)\right]ds\\ \nonumber
&+ \int_0^t\sum_{i, j=1}^v \mathbb E\left[\left(\hat H_{i} \hat x_{\pm} - \smat 0 \\ \hat c_i\srix\right)^T\hat \Sigma^{-1}\left(\hat H_{j} \hat x_{\pm} - \smat 0 \\ \hat c_j\srix\right)\right]k_{ij} ds.
\end{align}
Analogously to (\ref{Nkestimate2}), it holds that \begin{align}\label{blaNestimatg}
&\sum_{k=1}^m 2 \hat x_+^T(s)\hat\Sigma^{-1} \hat N_{k} \hat x_{\pm}(s) u_k(s) \\ \nonumber
&\leq  \hat x_+^T(s) \hat \Sigma^{-1} \hat x_+(s)  \left\|u^0(s)\right\|_{2}^2 +\sum_{k=1}^m  \hat x_{\pm}^T(s) \hat N_k^T\hat \Sigma^{-1} \hat N_{k} \hat x_{\pm}(s)\\ \nonumber
&\leq  \hat x_+^T(s) \hat \Sigma^{-1} \hat x_+(s)  \left\|u^0(s)\right\|_{2}^2 +\sum_{k=1}^m  \hat x_{\pm}^T(s) N_k^T \Sigma^{-1}  N_{k} \hat x_{\pm}(s).
\end{align}
Furthermore, we see that  
\begin{align*}
&\hat x_{\pm}^T (A^T\Sigma^{-1}+\Sigma^{-1} A)\hat x_{\pm}+4\hat x_{\pm}^T \Sigma^{-1} Bu = 2 \hat x_{\pm}^T \Sigma^{-1} (A\hat x_{\pm}+2Bu) \\
&= 2 \smat x_1 + x_r\\x_2 - h_1 \srix^T \hat\Sigma^{-1} (\hat A\hat x_{\pm}+2\hat B u) + 2 (\tilde x_3 - h_2)^T \Sigma_3^{-1}(\smat {A}_{31}&{A}_{32}& A_{33}\srix\hat x_{\pm}+2B_3u).
\end{align*}
Since  $\smat {A}_{31}&{A}_{32}& A_{33}\srix \smat x_{1}\\ x_{2}\\ \tilde x_{3}\srix= \smat {A}_{31}&{A}_{32}& A_{33}\srix \smat x_r\\ -h_1\\ -h_2\srix  =-B_3u$ by the definition of $\tilde x_3$
and the second line of (\ref{reltobeusef}), we obtain
$\smat {A}_{31}&{A}_{32}& A_{33}\srix\hat x_{\pm}=-2B_3 u$. Thus, \begin{align}\label{estimateforea22}
&\hat x_{\pm}^T (A^T\Sigma^{-1}+\Sigma^{-1} A)\hat x_{\pm}+4\hat x_{\pm}^T \Sigma^{-1} Bu \\ \nonumber 
&= 2 \hat x_+^T \hat\Sigma^{-1} (\hat A\hat x_{\pm}+2\hat B u) + 2\smat 0\\ -h_1\srix^T \hat \Sigma^{-1} (\hat A\hat x_{\pm}+2\hat B u).
\end{align}
Finally,  we see that
\begin{align}\label{relHterms}
 &\sum_{i, j=1}^v \left(\hat H_{i} \hat x_{\pm} - \smat 0 \\ \hat c_i\srix\right)^T\hat \Sigma^{-1}\left(\hat H_{j} \hat x_{\pm} - \smat 0 \\ \hat c_j\srix\right)k_{ij}\\ \nonumber
 &= \hat x_{\pm}^T\sum_{i, j=1}^v \hat H^T_{i}\hat \Sigma^{-1} \hat H_{j}k_{ij} \hat x_{\pm} - \sum_{i, j=1}^v \left(2 \hat H_{i} \hat x_{\pm} - \smat 0 \\ \hat c_i\srix\right)^T\hat \Sigma^{-1}\smat 0 \\ \hat c_j\srix k_{ij}\\\nonumber
 &\leq \hat x_{\mp}^T\sum_{i, j=1}^v H^T_{i} \Sigma^{-1} H_{j}k_{ij} \hat x_{\pm} - \sum_{i, j=1}^v \left(2 \hat H_{i} \hat x_{\pm} - \smat 0 \\ \hat c_i\srix\right)^T\hat \Sigma^{-1}\smat 0 \\ \hat c_j\srix k_{ij}
\end{align}
applying Lemma \ref{proppossemidef}. With (\ref{blaNestimatg}), (\ref{estimateforea22}) and (\ref{relHterms}) inequality (\ref{productruleappliedplus22}) becomes 
\begin{align}\nonumber
&\mathbb E\left[\hat x_+^T(t)\hat \Sigma^{-1} \hat x_+(t)\right]\\ \nonumber&\leq \mathbb E\int_0^t \hat x_{\pm}^T\left(A^T \Sigma^{-1}+\Sigma^{-1} A+\sum_{k=1}^m  N_k^T\Sigma^{-1} N_{k}
+\sum_{i, j=1}^v H^T_{i}\Sigma^{-1} H_{j}k_{ij}\right)\hat x_{\pm} ds\\ \label{insetobserveequation22}
&\quad -\mathbb E\int_0^t 2 \hat x_+^T\hat \Sigma^{-1} \smat 0 \\ \hat c_0\srix + \sum_{i, j=1}^v \left(2 \hat H_{i} \hat x_{\pm} - \smat 0 \\ \hat c_i\srix\right)^T \hat \Sigma^{-1}\smat 0 \\ \hat c_j\srix k_{ij}  ds\\
&\quad+ \mathbb E \int_0^t 2\smat 0 \\ h_1\srix^T \hat\Sigma^{-1} (\hat A\hat x_{\pm}+2 \hat Bu) ds + \mathbb E \int_0^t 4 \hat x_{\pm}^T\Sigma^{-1} Bu ds\nonumber \\
&\quad+\int_0^t \mathbb E\left[\hat x_+^T \hat \Sigma^{-1} \hat x_+\right]  \left\|u^0\right\|_{2}^2 ds. \nonumber
\end{align}
Similar to (\ref{ugrosserlyap}), we obtain
\begin{align*}
&4 \left\|u\right\|_{2}^2\geq \\ 
& \hat x_{\pm}^T\left(A^T \Sigma^{-1}+\Sigma^{-1} A+\sum_{k=1}^m  N_k^T\Sigma^{-1} N_{k}+\sum_{i, j=1}^v H^T_{i}\Sigma^{-1} H_{j}k_{ij}\right)\hat x_{\pm}+4 \hat x_{\pm}^T\Sigma^{-1} Bu.
\end{align*}
This leads to \begin{align}\nonumber
&\mathbb E\left[\hat x_+^T(t)\hat \Sigma^{-1} \hat x_+(t)\right]\\ \nonumber&\leq 4 \left\|u\right\|_{L^2_t}^2 +\int_0^t \mathbb E\left[\hat x_+^T \hat \Sigma^{-1} \hat x_+\right]  \left\|u^0\right\|_{2}^2 ds\\ \label{insetobserveezzquation22}
&\quad -\mathbb E\int_0^t 2 \hat x_+^T\hat \Sigma^{-1} \smat 0 \\ \hat c_0\srix + \sum_{i, j=1}^v \left(2 \hat H_{i} \hat x_{\pm} - \smat 0 \\ \hat c_i\srix\right)^T \hat \Sigma^{-1}\smat 0 \\ \hat c_j\srix k_{ij}  ds\\
&\quad+ \mathbb E \int_0^t 2\smat 0 \\ h_1\srix^T \hat\Sigma^{-1} (\hat A\hat x_{\pm}+2 \hat Bu) ds.\nonumber
\end{align}
In the following (\ref{insetobserveezzquation22}) is expressed by terms depending on $\Sigma_2$. We obtain $\hat x_+^T\hat \Sigma^{-1} \smat 0 \\ \hat c_0\srix=x_2^T\Sigma_2^{-1} \hat c_0$ exploiting
the partitions of $\hat x_+$ and $\hat \Sigma$. The terms depending on $\hat H_i$ become 
\begin{align} \nonumber
&-\sum_{i, j=1}^v\left(2 \hat H_{i} \hat x_{\pm} - \smat 0 \\ \hat c_i\srix\right)^T\hat \Sigma^{-1}\smat 0 \\ \hat c_j\srix k_{ij} 
=-\sum_{i, j=1}^v\left(2 \hat H_{i} \hat x_{\pm} - \smat 0 \\ \hat c_i\srix\right)^T\smat 0 \\ \Sigma_2^{-1} \hat c_j\srix k_{ij}\\  \nonumber
&=-\sum_{i, j=1}^v\left(2 \smat H_{i, 21} & H_{i, 22}& H_{i, 23}\srix (\smat x_1 \\ x_2\\\tilde x_3\srix + \smat x_r \\ -h_1\\-h_2\srix) -  \hat c_i\right)^T \Sigma_2^{-1} \hat c_j k_{ij}\\  \nonumber
&= -\sum_{i, j=1}^v\left(2 \smat H_{i, 21} & H_{i, 22}& H_{i, 23}\srix \smat x_1 \\ x_2\\\tilde x_3\srix +  \hat c_i\right)^T \Sigma_2^{-1} \hat c_j k_{ij}\\ \label{letzesh}
&\leq -\sum_{i, j=1}^v\left(2 \smat H_{i, 21} & H_{i, 22}& H_{i, 23}\srix \smat x_1 \\ x_2\\\tilde x_3\srix -  \hat c_i\right)^T \Sigma_2^{-1} \hat c_j k_{ij}
\end{align}
adding $\sum_{i, j=1}^v \hat c_i^T\Sigma_2^{-1} \hat c_j k_{ij}$ which is positive due to Lemma \ref{proppossemidef}.
Furthermore, using the first line of (\ref{reltobeusef}), it holds that \begin{align}\nonumber
&2\smat 0 \\ h_1\srix^T \hat \Sigma^{-1} (\hat A\hat x_{\pm}+2 \hat Bu) = 2\smat 0 & h_1^T\Sigma_2^{-1}\srix  (\hat A\hat x_{\pm}+2 \hat Bu)\\ \nonumber
&=  2 h_1^T\Sigma_2^{-1} (\smat A_{21} & A_{22}&A_{23}\srix (\smat x_1 \\ x_2\\\tilde x_3\srix + \smat x_r \\ -h_1\\-h_2\srix)+2 B_2u)   \\ \label{letztesa}
&=   2 h_1^T\Sigma_2^{-1} (\smat A_{21} & A_{22}&A_{23}\srix \smat x_1 \\ x_2\\\tilde x_3\srix + B_2u).                       
                           \end{align}
We insert (\ref{letzesh}) and (\ref{letztesa}) into (\ref{insetobserveezzquation22}) and obtain
\begin{align*}
\mathbb E\left[\hat x_+^T(t)\hat \Sigma^{-1} \hat x_+(t)\right]\leq 4 \left\|u\right\|_{L^2_t}^2 +\int_0^t \mathbb E\left[\hat x_+^T \hat \Sigma^{-1} \hat x_+\right]  \left\|u^0\right\|_{2}^2 ds-\hat \alpha_+(t).
\end{align*}
With Lemma \ref{gronwall}, analogously to (\ref{analoggronwallest}), we find \begin{align}\label{analoggronwallest22}
&\hat\alpha_+(t)+\int_0^t \hat \alpha_+(s)  \left\|u^0(s)\right\|_{2}^2 \exp\left(\int_s^t \left\|u^0(w)\right\|_{2}^2 dw\right) ds\\ \nonumber
&\leq 4\left\|u\right\|_{L^2_t}^2 \exp\left(\int_0^t \left\|u^0(s)\right\|_2^2ds\right).
\end{align}
The relations (\ref{boundmitalpmin}) and (\ref{analoggronwallest22}) yield the claim.
\end{proof}
\end{theorem}

\subsection{Proof of Theorem \ref{mainthmintro}}\label{proofmainthm}
We apply the results in Theorems \ref{mainthm} and \ref{mainthm2}. We remove the HSVs step by 
step and exploit the triangular inequality in order to bound the error between the outputs $y$ and $y_r$. We have\begin{align*}
  &\left(\mathbb E \left\|y-y_r\right\|^2_{L^2_T}\right)^{\frac{1}{2}}\\&\leq   
\left(\mathbb E\left\|y-y_{r_\nu}\right\|^2_{L^2_T}\right)^{\frac{1}{2}}+\left(\mathbb E\left\|y_{r_\nu}-y_{r_{\nu-1}}\right\|^2_{L^2_T}\right)^{\frac{1}{2}}+\ldots
+\left(\mathbb E\left\|y_{r_2}-y_{r}\right\|^2_{L^2_T}\right)^{\frac{1}{2}},        
  \end{align*}
where $y_{r_i}$ are the outputs of the ROMs with dimensions $r_i$ defined by $r_{i+1}=r_{i}+m(\tilde\sigma_{i})$ for $i=1, 2 \ldots, \nu-1$. Here, 
$m(\tilde\sigma_{i})$ denotes the multiplicity of $\tilde\sigma_{i}$ and $r_1=r$. In the reduction step from $y$ to $y_{r_\nu}$ only the smallest 
HSV $\tilde\sigma_\nu$ is removed from the system. Hence, by Theorem \ref{mainthm}, we have \begin{align*}                       
                  \left(\mathbb E \ \left\|y-y_{r_\nu}\right\|_{L^2_T}\right)^{\frac{1}{2}}\leq 2 \tilde\sigma_\nu \left\|u\right\|_{L^2_T}\exp\left(0.5 \left\|u^0\right\|_{L^2_T}^2\right).
\end{align*}
The ROMs of the outputs $y_{r_j}$ and $y_{r_{j-1}}$ are neighboring according to Section \ref{secneigeb}, i.e., only the HSV $\tilde\sigma_{r_{j-1}}$ is removed in the reduction step. 
By Theorem \ref{mainthm2}, we obtain\begin{align*}
  \left(\mathbb E \ \left\|y_{r_j}-y_{r_{j-1}}\right\|_{L^2_T}\right)^{\frac{1}{2}}\leq 2 \tilde\sigma_{r_{j-1}} \left\|u\right\|_{L^2_T} \exp\left(0.5 \left\|u^0\right\|_{L^2_T}^2\right)
  \end{align*}
for $j=2, \ldots, \nu $. This provides the claimed result.

\section{Conclusions}
In this paper, we investigated a large-scale stochastic bilinear system. In order to reduce the state space dimension, a model order reduction technique called singular perturbation approximation was extended to this setting. 
This method is based on Gramians proposed in \cite{redstochbil} that characterize how much a state contributes to the system dynamics. This choice of Gramians as well as the structure of the reduced system
is different than in \cite{hartmann}. With this modification, we provided a new $L^2$-error bound that can be used to point out the cases in which the reduced order model by singular perturbation approximation 
delivers a good approximation to the original model. This error bound is new even for deterministic bilinear systems.


\appendix

\section{Supporting Lemmas}

In this appendix, we state three important results and the corresponding references that we frequently use throughout this paper.
\begin{lemma}\label{lemstochdiff}
Let $a, b_1, \ldots, b_v$ be $\mathbb R^d$-valued processes, where $a$ is $\left(\mathcal F_t\right)_{t\geq 0}$-adapted and almost surely Lebesgue integrable and the functions $b_i$ 
are integrable with respect to the mean zero square integrable L\'evy process $M=(M_1, \ldots, M_v)^T$ with covariance matrix $K=\left(k_{ij}\right)_{i, j=1, \ldots, v}$. If the process $x$ is given by \begin{align*}
 dx(t)=a(t) dt+ \sum_{i=1}^v b_i(t)dM_i,                                                                  
                                                                   \end{align*}
then, we have \begin{align*}
 \frac{d}{dt}\mathbb E\left[x^T(t) x(t)\right]=2 \mathbb E\left[x^T(t) a(t)\right] + \sum_{i, j=1}^v \mathbb E\left[b_i^T(t) b_j(t)\right]k_{ij}.  
                                                                                                                                   \end{align*}
\begin{proof}
We refer to \cite[Lemma 5.2]{redmannspa2} for a proof of this lemma. 
\end{proof}
\end{lemma}
\begin{lemma}\label{proppossemidef}
Let $A_1, \ldots, A_v$ be $d_1\times d_2$ matrices and $K=(k_{ij})_{i, j=1, \ldots, v}$ be a positive semidefinite matrix, then
\begin{align*}\tilde K:=\sum_{i,j=1}^v A_i^T A_j k_{ij} \end{align*}
is also positive semidefinite.
 \begin{proof}
  The proof can be found in \cite[Proposition 5.3]{redmannspa2}.
 \end{proof}
\end{lemma}

\begin{lemma}[Gronwall lemma]\label{gronwall}
Let $T>0$, $z, \alpha: [0, T]\rightarrow \mathbb R$ be measurable bounded functions and $\beta: [0, T]\rightarrow \mathbb R$ be a nonnegative integrable function. 
If \begin{align*}
    z(t)\leq \alpha(t)+\int_0^t \beta(s) z(s) ds,
   \end{align*}
then it holds that \begin{align*}
    z(t)\leq \alpha(t)+\int_0^t \alpha(s)\beta(s) \exp\left(\int_s^t \beta(w)dw\right) ds
   \end{align*}
for all $t\in[0, T]$.
 \begin{proof}
  The result is shown as in \cite[Proposition 2.1]{gronwalllemma}. 
 \end{proof}
\end{lemma}

\bibliographystyle{plain}

\end{document}